\documentclass[11pt]{article}

\usepackage{amsmath}
\usepackage{amsfonts}
\usepackage{url}

\title{Simplicial blowups and discrete normal surfaces in the \texttt{GAP} package \texttt{simpcomp}}

\author{Felix Effenberger and Jonathan Spreer \\
Institut f\"ur Geometrie und Topologie \\
Universit\"at Stuttgart, Germany \\
\url{effenberger@mathematik.uni-stuttgart.de}, \\ \url{spreer@mathematik.uni-stuttgart.de}}

\date{\today}

\begin{document}

\maketitle

\begin{abstract}
	\texttt{simpcomp} is an extension to \texttt{GAP}, the well known system for computational discrete algebra. It allows the user to work with simplicial complexes. In the latest version, support for \emph{simplicial blowups} and  \emph{discrete normal surfaces} was added, both features unique to \texttt{simpcomp}. Furthermore, new functions for constructing certain infinite series of triangulations have been implemented and interfaces to other software packages have been improved to previous versions.
\end{abstract}

\section{Introduction}

\texttt{simpcomp} \cite{simpcompISSAC, simpcomp} is a package for working with simplicial complexes. Its aim is to provide the user with a broad spectrum of functionality regarding simplicial constructions and the calculation of properties of simplicial complexes.

Important goals during the development of \texttt{simpcomp} were interactivity, ease of use, completeness of documentation and ease of extensibility. The software allows the user to interactively construct simplicial complexes and to compute their properties in the \texttt{GAP} \cite{GAP4} or \texttt{SAGE} \cite{SAGE} shell. It is sought of as a tool for researchers to verify or disprove a conjecture one might have and to quickly do simplicial constructions using the computer. Furthermore, it makes use of \texttt{GAP}'s expertise in groups and group operations. For example, automorphism groups (cf.\ \cite{Bagchi10TrigCP3SymmCubeS2}) and fundamental groups of simplicial complexes can be computed and examined further within the \texttt{GAP} system.

%Apart from supplying a facet list, the user can as well construct simplicial complexes from a set of generators and a prescribed automorphism group -- the latter form being the common in which a complex is presented in a publication. This feature is to our knowledge unique to \texttt{simpcomp}. 
%
%\texttt{simpcomp} ships with an extensive library of known triangulations of manifolds and pseudomanifolds. This is the first time that they are easily accessible without having to look them up in the literature \cite{Kuehnel99CensusTight}, \cite{Casella01TrigK3MinNumVert}, or online \cite{Lutz08ManifoldPage} and thus allows the user to work with many different known triangulations without having to construct them first. In the current version 1.4.0, the library contains more than $7,000$ triangulations of manifolds and pseudomanifolds, including all vertex transitive triangulations of manifolds with few vertices of the classification in \cite{Lutz08ManifoldPage}.

With its development being started in March 2009, \texttt{simpcomp} still is a rather young project, but now already contains roughly 270 functions and its manual \cite{simpcomp} contains about 180 pages of documentation. At the ISSAC 2010 in Munich, \texttt{simpcomp} won the \emph{Best Software Presentation Award} by the Fachgruppe Computeralgebra. The package is currently under review by the GAP Council and is subject to acceptance as a \emph{shared package} of \texttt{GAP} that will be included in all future standard \texttt{GAP} installations.

Furthermore, the software system \texttt{SAGE} can be used to work with \texttt{simpcomp}. See \cite{simpcomp} for an interactive web demo of \texttt{simpcomp} based on \texttt{SAGE} notebooks. Connecting \texttt{simpcomp} more tightly with \texttt{SAGE} is planned for future releases.

The upcoming version 1.5 scheduled for May 2010 will be fully compliant to the standard \text{GAP} object mechanism and will have more advanced interfaces to other software packages.

\section{Why \texttt{simpcomp}}

\texttt{simpcomp} encapsulates all methods and properties of a simplicial complex in a new \texttt{GAP} object type (as an abstract data type). This way \texttt{simpcomp} can transparently cache properties already calculated, thus preventing unnecessary double calculations. This is mainly done by using the \texttt{GAP} native caching mechanism \cite{Breuer98GAP4TypeSystem}. It also takes care of the error-prone vertex labeling of a complex.

\texttt{simpcomp} is written entirely in the \texttt{GAP} scripting language. This has two implications:

(1) On the one hand, this limits the efficiency of the implementation as the scripted \text{GAP} code can never be as fast as native code. 

(2) On the other hand, this gives the user the possibility to see behind the scenes and to customize or alter \texttt{simpcomp} functions in an interactive way, profiting of all the mathematical and algebraic capabilities of the \texttt{GAP} scripting language. 

In the author's view, the advantages of (2) outweigh the drawbacks of (1). This was a major point when deciding on \texttt{simpcomp}'s design principles and also sets the software in contrast to other software packages like \texttt{polymake} \cite{Joswig00Polymake} that are heterogeneous, i.e.\ in which algorithms are implemented in various languages.

\section{\texttt{simpcomp} functionality}

\texttt{simpcomp}'s fundamental functions can be roughly divided into four groups: (i) functions constructing simplicial complexes, (ii) functions calculating properties of simplicial complexes, (iii) functions dealing with bistellar moves and (iv) functions concerning the library and the communication with other software packages -- for a full list of supported features see the documentation \cite{simpcomp} or use \texttt{GAP}'s built in interactive help system (all of \texttt{simpcomp}'s functions start with the prefix \texttt{SC}).

Concerning (i), complexes can be constructed by supplying a facet list or a set of generators together with a prescribed automorphism group -- the latter form being the common in which a complex is presented in a publication. This feature is to our knowledge unique to \texttt{simpcomp}. Furthermore, standard triangulations can by generated from scratch (simplex, cross polytope, cyclic polytope, stacked polytopes, etc.) and simplicial Cartesian products, connected sums, handle additions, etc. can be formed, enabling the user to obtain a wide variety of complexes with different properties. In (ii), basic properties and invariants of a simplicial complex like its dimension, the $f$-, $h$- or $g$-vector, Euler characteristic, (co-)homology groups, intersection form, Betti numbers, fundamental group, orientation, etc. can be computed. Concerning (iii): bistellar moves \cite{Pachner87KonstrMethKombHomeo, Bjoerner00SimplMnfBistellarFlips} allow to modify a given triangulation while leaving its PL homeomorphism type invariant (for an introduction to PL topology see \cite{Rourke72IntrPLTop, Kuehnel95TightPolySubm}). The concept has proven a powerful tool in combinatorial topology and can for example be used to reduce the vertex number of a given triangulation, to check if a simplicial complex is a manifold, to establish PL homeomorphisms between pairs of manifolds, to randomize complexes, to check whether a complex lies in a certain class of triangulations \cite{Effenberger11StackPolyTightTrigMnf}, and so on. Concerning (iv), there exist functions to save and load simplicial complexes to and from files (in an XML format) and to import and export complexes in various formats (e.g. from and to \texttt{polymake/TOPAZ} \cite{Joswig00Polymake}, \texttt{Macaulay2} \cite{M2}, \LaTeX, etc.). In addition, the internal library, currently containing more than $7,000$ triangulations, can be searched either using the name of a complex or a condition on the properties that it has to fulfill. The software also supports user libraries which can be used to organize own collections of triangulations produced with \texttt{simpcomp}.   

\section{New features in version 1.4}

\begin{itemize}
	\item Support for simplicial blowups: in algebraic geometry, blowups provide a useful way to study singularities of algebraic varieties \cite{Gompf}. The idea is to replace a point by all lines passing through that point. This concept is now available for combinatorial $4$-manifolds (cf.\ \cite{Spreer09CombPorpsOfK3}) and integrated into \texttt{simpcomp}. This functionality is to the authors' knowledge not provided by any other software package so far.
	\item Support for \emph{discrete normal surfaces} \cite{Kneser29ClosedSurfIn3Mflds, Haken61TheoNormFl, Spreer10NormSurfsCombSlic} and \emph{slicings}: slicings of combinatorial $d$-manifolds are (non-singular) $(d-1)$-dimensional level sets of polyhedral Morse functions. In dimension $3$, slicings are discrete normal surfaces. \texttt{simpcomp} supports discrete normal surfaces as a new object type and enables the user to generate and analyze slicings together with the corresponding Morse functions.
	\item New infinite series of highly symmetric triangulations: many highly symmetric triangulations occur as members of infinite series. Some of these series are well known and have been integrated into \texttt{simpcomp} already (simplex, cross polytope, cyclic polytope, etc.). Others were just recently found by the second author. \texttt{simpcomp} contains the first computer implementation of these series presented in \cite[Chapter 4]{Spreer10Diss}. 
	\item \texttt{homalg} interface: \texttt{simpcomp} now can use the \texttt{GAP} package \texttt{homalg} \cite{HomAlg} for its homology computations. This allows the computation of (co-)homology groups of simplicial complexes over arbitrary rings and fields, as well as the usage of all the functionality related to homological algebra that \texttt{homalg} provides.
\end{itemize}

\section{Roadmap: version 1.5 and beyond}

The current version of \texttt{simpcomp} is 1.4. On the roadmap for the upcoming versions are the following points.

\begin{itemize}
	\item Faster bistellar moves: currently, the algorithms to perform bistellar moves are implemented in the \texttt{GAP} scripting language. Since the performance of some of the algorithms implemented in \texttt{simpcomp} are mainly dependent on the running time of bistellar moves we plan to implement these functions in \texttt{C}. This should vastly speed up all calculations using bistallar moves. However, sticking to our design principles, the higher-level steering code related to bistellar moves will remain on the \texttt{GAP} side.
	\item We want to investigate how \texttt{simpcomp} can more closely interact with other software packages in the field (both \texttt{GAP} and non-\texttt{GAP}, e.g.\ \texttt{Macaulay2} and \texttt{SAGE}).
	\item There exists a combinatorial formula for calculating the Stiefel-Whitney class of combinatorial manifolds \cite{Banchoff79CombFormulaStiefelWhitney} due to Banchoff. We plan on including this, and possibly further invariants, into \texttt{simpcomp}.
	\item As a long term goal, we would like to provide the functionality to perform surgery on combinatorial $3$- and $4$-manifolds. This would be a step forward to constructing candidates for combinatorial manifolds with exotic PL structures as already done in the smooth setting by Akbulut \cite{Akbulut08CorksPlugsExoticStruc}.
\end{itemize}

\section{Examples}

This section contains a small demonstration of the capabilities of \texttt{simpcomp} in form of transcripts of the \texttt{GAP} shell for some example constructions. Most of the features presented below have been newly introduced in version 1.4.

\subsection{Normal surfaces in cyclic $4$-polytopes}

For $n \geq 3$, consider the cyclic $4$-polytope $C_4(2n)$ on $2n$ vertices with vertex labels $1$ to $n$. By Gale's evenness condition, neither the span of all odd nor the span of all even vertices in $C_4(2n)$ contains a triangle of $C_4(2n)$. Thus, given the combinatorial $3$-sphere $S=\partial C_4 (2n)$, a level set of a Morse function on $S$ separating the even from the odd vertices gives rise to a handle body decomposition of $S$ --- this is a discrete normal surface in the sense of \cite{Spreer10NormSurfsCombSlic}.

This construction can be done in \texttt{simpcomp} as follows. Note that we arbitrarily chose $n=5$ for demonstration purposes below. 
{\scriptsize
\begin{verbatim}
gap> LoadPackage("simpcomp");; #load the package
Loading simpcomp 1.4.0
by F.Effenberger and J.Spreer
http://www.igt.uni-stuttgart.de/LstDiffgeo/simpcomp
gap> c_4_10:=SCBdCyclicPolytope(4,10);
[SimplicialComplex

 Properties known: Chi, Dim, ... , TopologicalType, VertexLabels.

 Name="Bd(C_4(10))"
 Dim=3
 Chi=0
 F=[ 10, 45, 70, 35 ]
 HasBoundary=false
 Homology=[ [ 0, [ ] ], [ 0, [ ] ], [ 0, [ ] ], [ 1, [ ] ] ]
 IsConnected=true
 IsStronglyConnected=true
 TopologicalType="S^3"

/SimplicialComplex]
\end{verbatim}
}
\noindent
Above, we constructed the boundary of the cyclic $4$-polytope $\partial C_4 (10)$ on $10$ vertices. Note the properties of \texttt{c\_4\_10} already computed by \texttt{simpcomp}. We now look at the level set of a Morse function on \texttt{c\_4\_10} which separates even and odd vertices:
{\scriptsize
\begin{verbatim}
gap> sl:=SCSlicing(c,[[1,3,5,7,9],[2,4,6,8,10]]);
[NormalSurface

 Properties known: Chi, ConnectedComponents, ... , VertexLabels, Vertices.

 Name="slicing [ [ 1, 3, 5, 7, 9 ], [ 2, 4, 6, 8, 10 ] ] of Bd(C_4(10))"
 Dim=2
 Chi=-10
 F=[ 25, 70, 0, 35 ]
 IsConnected=true
 TopologicalType="(T^2)#6"

/NormalSurface]
\end{verbatim}
}
\noindent
The resulting polytopal complex on $25$ vertices is a discrete normal surface without triangles and with $35$ quadrilaterals. Topologically, it is the orientable surface with Euler characteristic $-10$, and thus homeomorphic to $(\mathbb{T}^2)^{\# 6}$. A triangulated version of this complex can be easily obtained as follows.
{\scriptsize
\begin{verbatim}
gap> trig:=SCNSTriangulation(sl);;
\end{verbatim}
}

\subsection{Combinatorial blowups of the Kummer variety $K^4$}

The $4$-dimensional abstract Kummer variety $K^4$ with 16 nodes leads to the $K3$ surface by resolving the $16$ singularities \cite{Spanier56HomKummerMnf}. Using \texttt{simpcomp}, this process could be carried out in a combinatorial setting for the first time, cf.\ \cite{Spreer09CombPorpsOfK3}. The first step of this so-called \emph{dilatation} or \emph{blowup process} can be done as follows.

We first load the singular $16$-vertex triangulation of $K^4$ due to K\"uhnel \cite{Kuhnel86MinTrigKummVar} from the library.
{\scriptsize
\begin{verbatim}
gap> SCLib.SearchByName("Kummer");
[ [ 7493, "4-dimensional Kummer variety (VT)" ] ]
gap> k4:=SCLib.Load(7493);        
[SimplicialComplex

 Properties known: AltshulerSteinberg, AutomorphismGroup, ...  
                   ... VertexLabels, Vertices.

 Name="4-dimensional Kummer variety (VT)"
 Dim=4
 AutomorphismGroupSize=1920
 AutomorphismGroupStructure="((C2 x C2 x C2 x C2) : A5) : C2"
 AutomorphismGroupTransitivity=1
 Chi=8
 F=[ 16, 120, 400, 480, 192 ]
 G=[ 10, 55, 60 ]
 H=[ 11, 66, 126, -19, 7 ]
 HasBoundary=false
 HasInterior=true
 Homology=[ [ 0, [ ] ], [ 0, [ ] ], [ 6, [ 2, 2, 2, 2, 2 ] ], [ 0, [ ] ], [ 1, [ ] ] ]
 IsCentrallySymmetric=false
 IsConnected=true
 IsEulerianManifold=true
 IsOrientable=true
 IsPM=true
 IsPure=true
 Neighborliness=2

/SimplicialComplex]
\end{verbatim}
}
\noindent
We now verify that the link of vertex $1$ in $K^4$ topologically is a real projective $3$-space. The ranks of its integral homology groups and its fundamental group are the following: 
{\scriptsize
\begin{verbatim}
gap> lk1:=k4.Link(1);;     
gap> lk1.Homology;
[ [ 0, [  ] ], [ 0, [ 2 ] ], [ 0, [  ] ], [ 1, [  ] ] ]
gap> pi:=lk1.FundamentalGroup;
<fp group with 61 generators>
gap> Size(pi);
2
\end{verbatim}
}
\noindent
We now verify that, as suspected, the complex is PL homeomorphic to the minimal $11$-vertex triangulation of $\mathbb{R}P^3$ from the library. This is done using a heuristic algorithm based on bistellar moves.
{\scriptsize
\begin{verbatim}
gap> SCLib.SearchByName("RP^3");
[ [ 45, "RP^3" ], [ 113, "RP^3=L(2,1) (VT)" ], ... ]
gap> minRP3:=SCLib.Load(45);;
gap> SCEquivalent(lk1,minRP3);
#I  SCReduceComplexEx: complexes are bistellarly equivalent.
true
\end{verbatim}
}
\noindent
Finally, we resolve the singularity of $K^4$ at vertex $1$ by a simplicial blowup.
{\scriptsize
\begin{verbatim}
gap> c:=SCBlowup(k4,1);;
#I  SCBlowup: checking if singularity is a combinatorial manifold...
#I  SCBlowup: ...true
#I  SCBlowup: checking type of singularity...
#I  SCReduceComplexEx: complexes are bistellarly equivalent.
#I  SCBlowup: ...ordinary double point (supported type).
#I  SCBlowup: starting blowup...
#I  SCBlowup: map boundaries...
#I  SCBlowup: boundaries not isomorphic, initializing bistellar moves...
#I  SCBlowup: found complex with smaller boundary: f = [ 15, 74, 118, 59 ].
...
#I  SCBlowup: found complex with smaller boundary: f = [ 11, 51, 80, 40 ].
#I  SCBlowup: found complex with isomorphic boundaries.
#I  SCBlowup: ...boundaries mapped succesfully.
#I  SCBlowup: build complex...
#I  SCBlowup: ...done.
#I  SCBlowup: ...blowup completed.
#I  SCBlowup: You may now want to reduce the complex via 'SCReduceComplex'.
\end{verbatim}
}
\noindent
Indeed, the second Betti number increased by $1$, again as expected.
{\scriptsize
\begin{verbatim}
gap> k4.Homology; 
[ [ 0, [  ] ], [ 0, [  ] ], [ 6, [ 2, 2, 2, 2, 2 ] ], [ 0, [  ] ], [ 1, [  ] ] ]
gap> c.Homology; 
[ [ 0, [  ] ], [ 0, [  ] ], [ 7, [ 2, 2, 2, 2 ] ], [ 0, [  ] ], [ 1, [  ] ] ] 
\end{verbatim}
}
\noindent
The resulting complex now only has $15$ singularities. By iterating this process $15$ more times, we obtain a combinatorial triangulation of the $K3$ surface with standard PL structure. 

\section{Acknowledgments}
The authors acknowledge support by the DFG: \texttt{simpcomp} has partly been developed within the DFG projects Ku 1203/5-2 and Ku 1203/5-3. 

{\scriptsize

\normalsize
\end{document}